# Identification of reference set and measurement of returns to scale in DEA: A least distance based framework

Mahmood Mehdiloozad [*]

*Department of Mathematics, College of Sciences, Shiraz University 71454, Shiraz, Iran*

Mohammad Bagher Ahmadi

*Department of Mathematics, College of Sciences, Shiraz University 71454, Shiraz, Iran*

[*] **Corresponding author**: M. Mehdiloozad
Ph.D. student
Department of Mathematics
College of Sciences
Shiraz University
Golestan Street | Adabiat Crossroad | Shiraz 71454 | Iran

E-mail: m.mehdiloozad@gmail.com
Tel.: +98.9127431689

# Identification of reference set and measurement of returns to scale in DEA: A least distance based framework

**Abstract** In data envelopment model (DEA), while either the most or least distance based frameworks can be implemented for targeting, the latter is often more relevant than the former from a managerial point of view due to easy attainability of the targets. To date, the two projection-dependent problems of reference set identification and returns to scale (RTS) measurement have been extensively discussed in DEA literature. To the best of our knowledge, nonetheless, there exists only one study which uses a closest projection for identifying reference set and accomplishes this task through a primal–dual linear programming based method. Motivated by this, we investigate the two aforementioned problems in a least distance based framework. First, we propose a *lexicographic multiple-objective programming* problem to find a unique closest projection for an inefficient decision making unit (DMU). Associated with the founded projection, we then identify all the possible closest reference DMUs by developing a linear programming model in the envelopment form. For an inefficient DMU, we next define *closest RTS* (CRTS) as the RTS of its least distance projection and measure the CRTS in two stages. Finally, we illustrate our proposed approach by a numerical example and compare the results with those found in the most distance based framework.

**Keywords** Data envelopment analysis; closest projection; lexicographic optimization; maximal closest reference set; closest returns to scale.

# 1 Introduction

Data envelopment analysis (DEA), introduced by Charnes, Cooper, and Rhodes (1978, 1979), has been proved to be an invaluable non-parametric technique for evaluating the relative efficiency of a set of homogenous decision making units (DMUs) with multiple inputs and multiple outputs. In DEA, the relative comparison is performed within a so-called production possibility set (PPS) that is empirically constructed based on observed data and some postulates. Relative to the constructed PPS, DEA identifies a DMU as either *efficient* or *inefficient* in the sense of Koopmans (1951). If a DMU found to be inefficient, then a unique or multiple projection(s) can be determined for it. Each projection suggests to where and by how much it should be improved in order to achieving full efficiency.

Though the determination of a projection can be made through either the *most* or *least* distance ways of targeting, the former was generally used in the traditional DEA literature because of computational ease. From a managerial point of view, however, the latter is often more relevant than the former since a closest projection is as similar as to the evaluated DMU and can be reached with less effort. That is, the smallest improvements in inputs and outputs of the evaluated DMU are required to reach the efficient targets (coordinates of the projection). For more details see, for example, Ando, Kai, Maeda, and Sekitani (2012), Aparicio and Pastor (2014), Aparicio, Ruiz, and Sirvent (2007), Baek and Lee (2009), Frei and Harker (1999), Fukuyama, Maeda, Sekitani, and Shi (2014), Gonzalez and Alvarez (2001), Jahanshahloo, Vakili, and Mirdehghan, (2012), Jahanshahloo, Vakili, and Zarepisheh, (2012), Lozano and Villa (2005), Pastor and Aparicio (2010) and Portela, Borges, and Thanassoulis (2003).

Two projection-dependent problems attracting considerable interest in the DEA literature are the identification of reference set and the measurement of returns to scale (RTS), which are focused in this study. The former problem concerns with the identification of a set of observed efficient DMUs, called *reference set*, for an inefficient DMU, against which this unit is directly compared for its efficiency improvement. This identification is crucially important from a managerial point of view, because it may not be meaningful in practice to introduce unobserved (virtual) targets as benchmarks for performance improvement. More theoretical and practical details on reference set identification can be found, for example, in Sueyoshi and Sekitani (2007a, 2007b), Krivonozhko, Førsund, and Lychev (2014), Roshdi, Mehdiloozad, and Margaritis (2014) and Mehdiloozad, Mirdehghan, Sahoo, and Roshdi (in press). The latter problem is related to determining the returns to scale (RTS) characterization of efficient points of the PPS. The economic concept of RTS was firstly introduced by Banker (1984) and Banker, Charnes, and Cooper (1984) into the DEA context and, since then, was extensively explored in the literature from both theoretical and practical aspects. See, for example, Banker and Thrall (1992), Tone (1996, 2005), Banker, Cooper, Seiford, Thrall and Zhu (2004), Sueyoshi and Sekitani (2007a, 2007b), Krivonozhko, Førsund, and Lychev (2012), Krivonozhko et al. (2014), and Mehdiloozad et al. (in press), among others.

Based on the similarity between an inefficient DMU and its closest projection(s), the use of a closest projection in either of the two above-mentioned problems seems to be more meaningful than that of a farthest one. To the best of our knowledge, however, none of the existing RTS measurement methods has employed the least distance way of target setting. Moreover, among the previous research studies conducted on the reference set identification, the recent work by Roshdi et al. (2014) is the only one which uses a least distance projection for identifying *maximal closest reference set* (MCRS) – a concept defined as the set of all possible reference DMUs associated with a given closest projection. To identify a closest efficient projection, Roshdi et al. (2014) exploited the single-stage approach. Moreover, they employed the primal–dual based approach of Sueyoshi and Sekitani (2007a) for finding the MCRS.

The contribution of this paper is three-fold. First, modifying the innovative approach of Aparicio et al. (2007), we develop a *lexicographic multiple-objective programming (LMOP)*[1] problem to generate a *unique* least distance projection for an inefficient DMU. The LMOP problem minimizes the distance between the under evaluation DMU and the efficient frontier in a hierarchical manner and based on the priority rankings[2] of minimizing input excesses and output shortfalls. Note that while the approach of Aparicio et al. (2007) uses a single mixed integer linear programming (MILP) problem to find a closest projection, our lexicographic approach, on the contrary, requires solving a series of MILP problems in the number of inputs and outputs. Nonetheless, our approach is advantageous since its lexicographic nature guarantees the uniqueness of the obtained closest projection, which is our main motivation in developing our approach. As we illustrate via an example, this uniqueness is of crucial importance because the occurrence of multiple closest projections for an inefficient DMU may results in different types of RTSs and different MCRSs, each associated with a closest projection. In addition, the multiplicity issue in the least distance based framework cannot be deal effectively as in the most distance based one. This is due to the fact that multiple closest projections may not be located on the intersection of some efficient facets and may not constitute a face of the PPS, accordingly.

---

[1] Lozano and Villa (2009) proposed a target setting DEA approach based on the lexicographic multi-objective optimization.

Next, we define the notion of *unary closest reference set* (UCRS) as the set of efficient DMUs that are active in a specific convex combination which produces the projection obtained via the LMOP problem. We also redefine the notion of *maximal closest reference set* (MCRS) as the union of all the UCRSs. Based on the work of Mehdiloozad et al. (in press), we then propose a linear programming (LP) model in the envelopment form to identify the MCRS as our second contribution. The computational efficiency of our approach is higher than that of the primal–dual method of Roshdi et al (2014). This is because our proposed LP problem is formulated based on the primal (envelopment) form that is computationally more efficient than the dual (multiplier) form (Cooper, Seiford, and Tone, 2007). Furthermore, since our proposed LP problem contains several upper-bounded variables, the computational efficiency of our method can be enhanced by using the simplex algorithm[3] adopted for solving the LP problems with upper-bounded variables.

---

[2] To elicit preference information from the decision maker, the analytic hierarchy process (AHP) of Saaty (1980) can be used.

As pointed out by Banker et al. (2004), the RTS generally has an unambiguous meaning for an efficient DMU. Nonetheless, an efficient projection of an inefficient DMU is used in order to estimate its RTS. To the best of our knowledge, all the existing RTS measurement methods are farthest projection based. Therefore, the relevance of using a closest projection motivated us to introduce the notion of *closest RTS (CRTS)* as the third contribution of this paper. The CRTS is defined as the RTS of our founded unique closest projection of the evaluated inefficient DMU. We then use the method of Banker et al. (2004) to measure the RTS at this projection.

The remainder of the paper unfolds as follows. Section 2 recalls some preliminaries concerning the measurement of RTS. Section 3, first, clarifies the notion of the CRTS by the aid of a motivating example; second, it formulates an LMOP problem to find out a unique closest projection on the efficient frontier, and proposes an LP problem to identify the MCRS; finally, it introduces the notion of CRTS and develops a method for its determination. Section 4 presents a numerical example, and Section 5 concludes.

## 2 Preliminaries

Throughout this paper, we deal with *n* observed DMUs; each uses *m* inputs to produce *s* outputs. Let $\mathbf{x}_j = (x_{1j},...,x_{mj})^T \in \mathbb{R}^m_{\geq 0}$ and $\mathbf{y}_j = (y_{1j},...,y_{sj})^T \in \mathbb{R}^s_{\geq 0}$ represent the input and output vectors, respectively, for the *j*th DMU where $j \in J = \{1,...,n\}$. The superscript *T* stands for a vector transpose, and $o \in J$ represents the index of the under evaluation DMU.

DMU$_o$ is assessed with respect to the so-called *production possibility set (PPS)* defined as

$$T = \{(\mathbf{x}, \mathbf{y}) \in \mathbb{R}^m_{\geq 0} \times \mathbb{R}^s_{\geq 0} \mid \mathbf{x} \text{ can produce } \mathbf{y}\}. \tag{1}$$

---

[3] The simplex algorithm for bounded variables was published by Dantzig (1955) and was independently developed by Charnes and Lemke (1954). The use of this algorithm is much more efficient than the ordinary simplex algorithm for solving the LP problem with upper-bounded variables (Winston 2003).

Under variable returns to scale (VRS) framework, the non-parametric DEA-based representation of $T$ is set up as follows (Banker et al., 1984):

$$T_V = \left\{(\mathbf{x}, \mathbf{y}) \mid \sum_{j \in J} \lambda_j \mathbf{x}_j \leq \mathbf{x}, \ \sum_{j \in J} \lambda_j \mathbf{y}_j \geq \mathbf{y}, \ \sum_{j \in J} \lambda_j = 1, \ \lambda_j \geq 0, \forall j \in J \right\}. \tag{2}$$

In reference to $T_V$, the *envelopment* (primal) form of the BCC model is as

$$\begin{aligned}
\min \quad & \theta - \varepsilon \left( \sum_{i=1}^m s_i + \sum_{r=1}^s s_{m+r} \right) \\
\text{s.t.} \quad & \sum_{j \in J} \lambda_j x_{ij} + s_i = \theta x_{io}, \quad i = 1, \ldots, m, \\
& \sum_{j \in J} \lambda_j y_{rj} - s_{m+r} = y_{ro}, \quad r = 1, \ldots, s, \\
& \sum_{j \in J} \lambda_j = 1, \\
& \lambda_j \geq 0, \forall j \in J, \ s_i, s_{m+r} \geq 0, \forall i, r,
\end{aligned} \tag{3}$$

where $\varepsilon$ is a non-Archimedean infinitesimal. The dual of model (3), called the *multiplier* form, is as

$$\begin{aligned}
\max \quad & \sum_{r=1}^s w_{m+r} y_{ro} - w_0 \\
\text{s.t.} \quad & \sum_{i=1}^m w_i x_{io} = 1, \\
& \sum_{r=1}^s w_{m+r} y_{rj} - \sum_{i=1}^m w_i x_{ij} - w_0 \leq 0, \quad \forall j \in J, \\
& w_i \geq \varepsilon, \ w_{m+r} \geq \varepsilon, \ \forall i, r, \ w_0 : \text{free in sign}.
\end{aligned} \tag{4}$$

Let $(\theta^*, \boldsymbol{\lambda}^*, \mathbf{s}^*)$ be an optimal solution to model (3). Then, DMU$_o$ is called *efficient* if and only if $\theta^* = 1$ and $\mathbf{s}^* = \mathbf{0}_{m+s}$. Equivalently, DMU$_o$ is efficient if and only if $\sum_{r=1}^s w^*_{m+r} y_{ro} - w^*_0 = 1$ holds for an optimal solution $(\mathbf{w}^*, w^*_0)$ to model (4). If DMU$_o$ is inefficient, then its projection on the efficient frontier is given by

$$(\bar{\mathbf{x}}_o, \bar{\mathbf{y}}_o) := \sum_{j \in J} \lambda_j^* (\mathbf{x}_j, \mathbf{y}_j). \tag{5}$$

In model (4), the free variable $w_0$ corresponds to the convexity constraint $\sum_{j \in J} \lambda_j = 1$ in (3) and its optimal value is the intercept of the supporting hyperplane $\sum_{r=1}^{s} w_{m+r}^* y_r - \sum_{i=1}^{m} w_i^* x_i - w_0^* = 0$ binding at DMU$_o$. If DMU$_o$ is efficient, then the following theorem indicates that its RTS can be determined by examining the sign of $w_0$.

**Theorem 2.1** Let DMU$_o$ be efficient and let $\bar{w}_0$ and $\underline{w}_0$ be the *upper* and *lower* bound of $w_0$, respectively, which can be obtained by solving the following model:

$$\begin{aligned}
\bar{w}_0 (\underline{w}_0) = \max(\min) \quad & w_0 \\
\text{s.t.} \quad & \sum_{i=1}^{m} w_i x_{io} = 1, \\
& \sum_{r=1}^{s} w_{m+r} y_{rj} - \sum_{i=1}^{m} w_i x_{ij} - w_0 \leq 0, \quad \forall j \in J, \\
& \sum_{r=1}^{s} w_{m+r} y_{ro} - \sum_{i=1}^{m} w_i x_{io} - w_0 = 0, \\
& w_i \geq 0, \ w_{m+r} \geq 0, \ \forall i, \ r, \ w_0 : \text{free in sign}.
\end{aligned} \tag{6}$$

Then,

- IRS prevail at DMU$_o$ if $\bar{w}_0 < 0$.
- CRS prevail at DMU$_o$ if $\underline{w}_0 \leq 0 \leq \bar{w}_0$.
- DRS prevail at DMU$_o$ if $\underline{w}_0 > 0$.

Based on this theorem, the identification of RTS can be made in two stages. Note that if DMU$_o$ is inefficient, then its RTS is defined as the RTS of its projection given in (5).

# 3 Our proposed approach

## 3.1 Motivating example

Consider a data set consisting of four hypothetical DMUs labeled as A, B, C, and D, where each unit uses one input to produce one output. The input–output data are exhibited in Table 3.1.1 and are depicted in Fig. 3.1.1.

**Table 3.1.1** Data for four DMUs

|   | A | B | C | D |
|---|---|---|---|---|
| x | 2 | 3 | 6 | 4 |
| y | 2 | 5 | 6 | 4 |

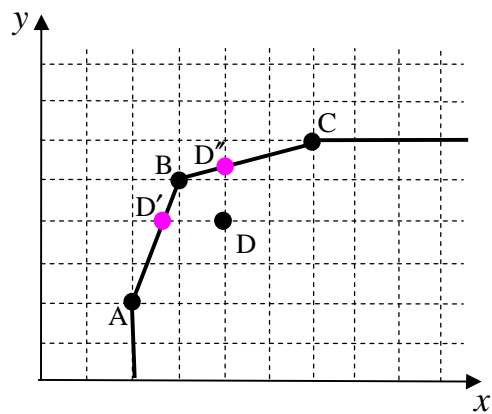

**Fig. 3.1.1** Graphical representation of the data in Table 3.1.1

As can be seen in Fig. 3.1.1, all DMUs are efficient except D. It can be easily shown that the efficient units A, B and C exhibit IRS, CRS and DRS. Furthermore, IRS and DRS respectively prevail for the relative interior points of the line segments AB and BC forming the efficient frontier.

To determine the RTS of the inefficient unit D, we need to project it onto the efficient frontier. In this example, all the points on the line segments D'B and BD" could potentially used as projections for D. If the range-adjusted model (RAM) of Cooper, Park, and Pastor (1999) is employed, then B is the unique furthest projection. By using the mRAM model of Aparicio et al. (2007), in contrast, $D'$ and $D''$ are identified as closest projections for D. Therefore, it is easier for D to attain to either $D'$ or $D''$ rather than B.

With regard to the similarity between D and its closest projections, it is more meaningful to use either $D'$ or $D''$ for identifying the reference DMUs and for estimating the RTS of D. Based on closest projections D' and D", the RTSs estimated for D are IRS and DRS and the corresponding reference sets are $\{A,B\}$ and $\{B,C\}$, respectively.

Although there are some DEA models that provide closest targets, the problem underlying the use of these models is the occurrence of multiple projections; because, as illustrated above, such an occurrence may lead to different reference sets and different RTSs. This verifies the need to use a unique closest projection for determining reference set and for measuring RTS. Note that the occurrence of multiple closest projections cannot be dealt as that of furthest ones. This is because the closest projection set may not be convex and the multiple closest projections may not be located on the intersection of efficient facet (see $D'$ and $D''$ in Fig 3.1.1).

In view of the above discussion, we make a modification on the approach of Aparicio et al. (2007) to obtain a unique closest projection. In this regard, we develop a lexicographic multiple-objective programming (LMOP) problem in the immediately following subsection that enable us to uniquely project an inefficient unit onto the efficient frontier. The LMOP problem accomplishes the task by minimizing the sum of all input excesses and output shortfalls based on the priority rankings provided by the decision maker, but not by minimizing the weighted sum of them as in the mRAM model.

## 3.2 Finding a unique closest projection

As illustrated in the previous subsection, the closest projection(s) provided by the mRAM model may not be uniquely determined for an inefficient DMU. To deal with this issue, we define $\Omega_o$ as the set of all input–output slack vectors required to project DMU$_o$ onto the efficient frontier or, formally, as

$$\Omega_o := \left\{ \mathbf{s} \,\middle|\, (\mathbf{x}_o, \mathbf{y}_o) - \mathbf{s} \in \partial^S(T_V) \right\}, \tag{7}$$

where $\partial^S(\cdot)$ represents the efficient frontier of $T_V$.

To compare the slack vectors in $\Omega_o$, we consider the lexicographic order[4]. That is, first, we arrange the input and output slacks in a hierarchical manner and rank them according to their priorities. Next, we minimize the distance between DMU$_o$ and $\partial^S(T_V)$ by finding a lexicographic solution of the following LMOP problem:

$$\operatorname*{lexmin}_{\mathbf{s}\in\Omega_o} \mathbf{s}. \tag{8}$$

By (7) and the theorem in Aparicio et al. (2007), we rewrite problem (8) as

$$\begin{aligned}
\operatorname{lexmin} \quad & \mathbf{s} \\
\text{s.t.} \quad & \sum_{j\in J_E} \lambda_j x_{ij} = x_{io} - s_i, \quad i=1,\ldots,m, \\
& \sum_{j\in J_E} \lambda_j y_{rj} = y_{ro} + s_{m+r}, \quad r=1,\ldots,s, \\
& \sum_{j\in J_E} \lambda_j = 1, \\
& \sum_{r=1}^{s} w_{m+r} y_{rj} - \sum_{i=1}^{m} w_i x_{ij} - w_0 + d_j = 0, \ \forall j\in J_E, \\
& 0 \le \lambda_j \le M(1-I_j), \ \forall j\in J_E, \\
& d_j - M I_j \le 0, \ I_j \in \{0,1\}, \ \forall j\in J_E, \\
& w_i \ge 1, \ s_i \ge 0, \ w_{m+r} \ge 1, \ s_{m+r} \ge 0, \ \forall i,r,
\end{aligned} \tag{9}$$

where $J_E$ denotes the index set of efficient DMUs, and "lexmin" represents lexicographical minimization.

Let $(\boldsymbol{\lambda}^*, \mathbf{s}^*, \mathbf{w}^*, w_0^*, \mathbf{d}^*, \mathbf{I}^*)$ be a lexicographic solution to model (9). Then, the unique closest projection for DMU$_o$, $P_o^C$, is obtained as

$$P_o^C = (\hat{\mathbf{x}}_o, \hat{\mathbf{y}}_o) = \sum_{j\in J} \lambda_j^* (\mathbf{x}_j, \mathbf{y}_j) = (\mathbf{x}_o, \mathbf{y}_o) - \mathbf{s}^*. \tag{10}$$

---

[4] In Ehrgott (2005), see Table 1.2 for a definition of the lexicographic order and Section 5.1 for details on lexicographic optimization.

This projection represents lexicographically less demanding levels of operation for the inputs and outputs of DMU$_o$ to perform efficiently. Note that despite model (9) may have multiple optimal solutions, the optimal vector $\mathbf{s}^*$ is the same for all of these solutions. This indicates that $P_o^C$ is determined uniquely by all of the possible optimal solutions.

### 3.3 Identification of maximal closest reference set

From Mehdiloozad et al. (in press) and Roshdi et al. (2014), we present the following definition.

**Definition 3.3.1** Let $\left(\boldsymbol{\lambda}^*, \mathbf{s}^*, \mathbf{w}^*, w_0^*, \mathbf{d}^*, \mathbf{I}^*\right)$ be a lexicographic solution to model (9). Then, the set of DMUs with positive $\lambda_j^*$ is defined as the *unary closest reference set (UCRS)* for DMU$_o$ and is denoted by $R_o^C$ as

$$R_o^C := \left\{ \mathrm{DMU}_j \mid \lambda_j^* > 0 \right\}. \tag{11}$$

Since $P_o^C$ may be expressed as multiple convex combinations of closest reference DMUs, multiple values may take place for the vector $\boldsymbol{\lambda}$ which can lead to the occurrence of multiple UCRSs. To deal with an occurrence of multiple UCRSs, we define a unique closest reference set containing all the possible UCRSs.

**Definition 3.2.2** The set of *all* the UCRSs is defined as the *maximal closest reference set (MCRS)* for DMU$_o$ and is denoted by $R_{oP}^{MC}$ as

$$R_o^{MC} := \left\{ \mathrm{DMU}_j \mid \lambda_j^* > 0 \text{ in some lexicographic solution of (9)} \right\}. \tag{12}$$

Note that there exists a face of minimum dimension, called *minimum face*, which contains $P_o^C$. Moreover, $P_o^C$ lies in the relative interior of this face; otherwise, there would be a face of a dimension less than that of the minimum face that contains $P_o^C$. Therefore, we can conclude that Lemma 3.1.1, Theorem 3.1.1 and its two corollaries in Mehdiloozad et al. (in press) similarly hold true for the MCRS. Furthermore, analogous to Theorem 3.1.2 in Mehdiloozad et al. (in press), the following result can be derived.

**Theorem 3.3.1** Let $\boldsymbol{\lambda}^{\max}$ be a solution to the following system of constraints for which the number of positive components is maximum:

$$\sum_{j \in J_E} \lambda_j \mathbf{x}_j = \hat{\mathbf{x}}_o,$$
$$\sum_{j \in J_E} \lambda_j \mathbf{y}_j = \hat{\mathbf{y}}_o, \qquad (13)$$
$$\sum_{j \in J_E} \lambda_j = 1, \ \lambda_j \geq 0, \ \forall j \in J_E,$$

where $(\hat{\mathbf{x}}_o, \hat{\mathbf{y}}_o)$ is the closest projection given in (10). Then,

$$R_o^{MC} = \left\{ \mathrm{DMU}_j \,\middle|\, \lambda_j^{\max} > 0 \right\}. \qquad (14)$$

Based on the above theorem and Theorem 3.2.1 of Mehdiloozad et al. (in press), we now propose the following LP problem in order to identify $R_o^{MC}$:

$$\begin{aligned}
\max \ & \sum_{k=1}^{t+1} \alpha_{j_k} \\
\text{s.t.} \ & \sum_{k=1}^{t} \left( \alpha_{j_k} + \beta_{j_k} \right) \mathbf{x}_j - \left( \alpha_{j_{t+1}} + \beta_{j_{t+1}} \right) \hat{\mathbf{x}}_o = \mathbf{0}_m, \\
& \sum_{k=1}^{t} \left( \alpha_{j_k} + \beta_{j_k} \right) \mathbf{y}_j - \left( \alpha_{j_{t+1}} + \beta_{j_{t+1}} \right) \hat{\mathbf{y}}_o = \mathbf{0}_s, \qquad (15)\\
& \sum_{k=1}^{t} \left( \alpha_{j_k} + \beta_{j_k} \right) - \left( \alpha_{j_{t+1}} + \beta_{j_{t+1}} \right) = 0, \\
& 0 \leq \alpha_{j_k} \leq 1, \ \beta_{j_k} \geq 0, \ \forall k = 1,\ldots,t+1,
\end{aligned}$$

where $t$ denotes the cardinality of $J_E$, i.e., $J_E = \{j_1, \ldots, j_t\}$.

Let $(\boldsymbol{\alpha}^*, \boldsymbol{\beta}^*)$ be an optimal solution to model (15). Then, by Theorem 3.2.1 in Mehdiloozad et al. (in press), we have

$$\lambda_{j_k}^{\max} = \frac{1}{\alpha_{j_{t+1}}^* + \beta_{j_{t+1}}^*} \left( \alpha_{j_k}^* + \beta_{j_k}^* \right), \quad k = 1, \ldots, t. \tag{16}$$

In accordance with Theorem 3.3.1, $R_o^{MC}$ can be hence identified via model (15), and $P_o^C$ can be expressed as a strict convex combination of the reference DMUs in this set, accordingly.

Note that the computational efficiency of our proposed approach is higher than that of the primal–dual one proposed by Roshdi et al. (2014), since ours is developed based on the primal (envelopment) form which is computationally more efficient than the dual (multiplier) form (Cooper et al., 2007). Moreover, since our proposed LP problem contains several upper-bounded variables, the computational efficiency of our approach can be enhanced by using the simplex algorithm adopted for solving the LP problems with upper-bounded variables.

### 3.4 Measurement of closest returns to scale (CRTS)

In the conventional RTS measurement methods, the RTS of an inefficient DMU is commonly defined as that of a furthest projection of this unit. However, from a target setting perspective, the more the projection is near to an assessed inefficient DMU, the less levels of operation for inputs and outputs of this DMU are needed to make it efficient. In this sense, the projection attained by the smallest modifications in inputs and outputs of the assessed unit is as much similar as to the assessed DMU. Based on this fact, we present the following definition.

**Definition 3.4.1** Let $P_o^C$ be the unique closest projection of DMU$_o$ as defined in (10). Then, we define the *closest RTS (CRTS)* for DMU$_o$ as the RTS of $P_o^C$.

We now proceed to develop a method for determining the CRTS. It is known that both the type and magnitude of RTS for $P_o^C$ can be determined through the position(s) of the hyperplane(s) supporting the PPS at this projection. The supporting hyperplane(s) passes/pass through the CMRS associated with $P_o^C$ and can be mathematically characterized by this CMRS. Nonetheless, the uniqueness of the characterized supporting hyperplane(s) cannot be guaranteed since the face of minimum dimension that contains this projection may not be an efficient facet of dimension $m+s-1$ in the input–output space, i.e., the minimum face may not be a 'Full Dimensional Efficient Facet' (Olesen and Petersen, 1996, 2003). To deal with the occurrence of multiple supporting hyperplanes, an interesting and effective method is that of Banker et al. (2004) as described in Section 2. Using this method, we measure the CRTS at $P_o^C$.

## 4 Numerical example

In this section, we apply our proposed approach to the data of eight hypothetical DMUs with one input and one output. The data taken from Mehdiloozad et al. (in press) is exhibited in Table 4.1 and its generated frontier is depicted in Fig. 4.1.

**Table 4.1** Input–output data

|        | $DMU_1$ | $DMU_2$ | $DMU_3$ | $DMU_4$ | $DMU_5$ | $DMU_6$ | $DMU_7$ | $DMU_8$ |
|--------|---------|---------|---------|---------|---------|---------|---------|---------|
| **Input**  | 1 | 2 | 3 | 5 | 8 | 2 | 3 | 6 |
| **Output** | 2 | 5 | 6 | 8 | 8 | 1 | 3 | 4 |

**Source:** Mehdiloozad et al. (in press)

As can be seen in Fig. 4.1, $DMU_1$, $DMU_2$, $DMU_3$ and $DMU_4$ are efficient, and the remaining units are inefficient. Moreover, $J_E = \{1,2,3,4\}$ and the efficient frontier consists of the two line segments connecting $DMU_1$ to $DMU_2$ and $DMU_2$ to $DMU_4$. It can be easily shown that the efficient units $DMU_1$ and $DMU_2$ exhibit IRS and CRS, respectively, and $DMU_3$ and $DMU_4$ both exhibit DRS.

To determine the MCRS and the CRTS for each inefficient unit, we first project it onto the efficient frontier via the LMOP model. To employ this model, suppose that the given priority ranking of input and output is $2$ and $1$. That is, first, the output slack is minimized; then, by fixing the output slack at its optimum value, the input slack is minimized.

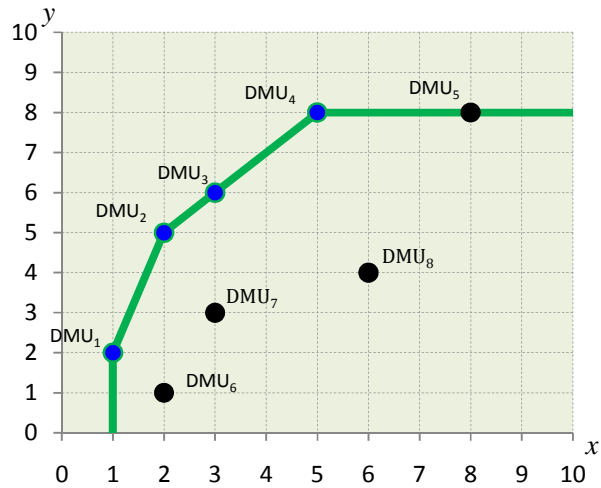

**Fig 4.1** The production frontier

Table 4.2 reports the results for the inefficient DMUs. While $DMU_4$ and $DMU_1$ are, respectively, the closest observed projections for the inefficient units $DMU_5$ and $DMU_6$, the closest projections of $DMU_7$ and $DMU_8$ are the unobserved points (1.3333,3) and (1.6667,4), respectively.

**Table 4.2** The MCRSs and the CRTSs for inefficient units

|  |  |  | $DMU_5$ | $DMU_6$ | $DMU_7$ | $DMU_8$ |
|---|---|---|---|---|---|---|
| **LMOP Model (9)** |  |  | $DMU_4$ | $DMU_1$ | (1.333,3) | (1.667,4) |
| **Model (15)** | Ref. Weights | $\lambda'_1$ |  | 1 | 0.6667 | 0.3333 |
|  |  | $\lambda'_2$ |  |  | 0.3333 | 0.6667 |
|  |  | $\lambda'_3$ |  |  |  |  |
|  |  | $\lambda'_4$ | 1 |  |  |  |
| **CRTS measurement** |  | $u_o^*$ | 0.600 | -1 | -0.333 | -0.333 |
|  |  | $u_o^+ / u_o^-$ | 0.600 | -0.3333 | -0.333 | -0.333 |
|  |  | **CRTS** | DRS | IRS | IRS | IRS |

**CRS:** Constant Returns to Scale; **IRS:** Increasing Returns to Scale; **DRS:** Decreasing Returns to Scale.

Since $DMU_4$ and $DMU_1$ are both extreme-efficient units, so the MCRSs of $DMU_5$ and $DMU_6$ are their respective closest projections. However, using model (15) for the points (1.3333, 3) and (1.6667, 4) indicates that the MCRS for both of the units $DMU_7$ and $DMU_8$ is $\{DMU_1, DMU_2\}$. Hence, as we expect, the minimum face containing each of the points (1.3333, 3) and (1.6667, 4) is the line segment connecting $DMU_1$ to $DMU_2$. However, as recorded in Table 2 in Mehdiloozad et al. (2014), the furthest reference set of $DMU_6$ is $DMU_2$. Furthermore, the furthest reference set of $DMU_7$ and $DMU_8$ consists of $DMU_2$, $DMU_3$ and $DMU_4$.

We now turn to estimate the CRTS of each inefficient DMU based upon the RTS of its closest projection. To do so, we use the two–stage method of Banker et al. (2004), and examine the sign of the intercept of the supporting hyperplane passing through the closest projection of the evaluated unit. The results are summarized in the last three lines of Table 4.2. While the inefficient units $DMU_6$, $DMU_7$ and $DMU_8$ have increasing CRTS, the CRTS of $DMU_5$ is decreasing. However, as given in Table 2 in Mehdiloozad et al. (2014), if the RTSs for the inefficient units are measured based on their furthest RAM-projections, then $DMU_5$, $DMU_7$ and $DMU_8$ exhibit decreasing RTS and constant RTS prevails at $DMU_6$.

## 5 Concluding remarks

Two frameworks can be applied for targeting in DEA. The first and commonly used one is the most distance based framework that provides furthest projections for the under evaluation DMU. From a managerial point of view, this framework may be not appropriate when the decision maker wishes that the resulting targets to be as similar as to inputs and outputs of the evaluated DMU. The second one that is much relevant in such situations is the least distance based framework that produces closest projections. The desirability of this framework over the first one is due to the fact that the targets associated with closest projections are easily attainable.

In view of the above discussion, the current study was concerned with establishing a least distance based framework to discuss the two projection-dependent issues of reference set identification and RTS measurement. In this regard, first, we exploited the interesting approach of Aparicio et al. (2007) and developed an LMOP problem to provide a unique least distance projection for an evaluated inefficient DMU. On comparison between our approach and that of Aparicio et al. (2007), one may argue that while the former requires solving a series of MILP problems, the latter uses a single MILP problem to find a closest projection. Nonetheless, the former is particularly advantageous to the latter in terms of guaranteeing the uniqueness of the obtained projection. This is because the occurrence of multiple closest projections in the approach of Aparicio et al. (2007) may result in different types of RTS and different MCRSs, which is not meaningful.

Second, based on the works of Mehdiloozad et al. (2014) and Roshdi et al. (2014), we redefined two notions: i) UCRS: the set of efficient DMUs that are active in a specific convex combination generating the unique projection obtained via the LMOP problem, and ii) MCRS: the union of all UCRSs associated with this projection. Then, we proposed an LP problem to identify the MCRS. Since this LP problem involves some upper-bounded variables, the computational efficiency of our approach would be improved using the simplex algorithm adopted for solving the LP problems with upper-bounded variables. Moreover, as the developed problem is a primal-based LP, it has less number of constraints that improves the computational efficiency compared to the existing primal–dual based approach.

Third, we introduced the notion of CRTS that is defined as the RTS of the evaluated DMU's closest projection. To measure the CRTS, we used the method of Banker et al. (2004), which can effectively deal with the occurrence of multiple supporting hyperplanes due to non-full dimensionality of the minimum face.

To ensure the uniqueness of the obtained closest projection, we incorporated the decision maker's preferences through the lexicographic optimization method, which requires a priori articulation of preferences. Hence, our proposed approach can be used when the decision maker is able to state his preferences a priori. Since interactive methods enable the decision maker to incorporate his/her preferences iteratively, an interesting topic for future research would be to develop an interactive method for finding a unique closest efficient target.

Krivonozhko, V. E., Førsund, F. R., & Lychev, A. V. (2012). Returns-to-scale properties in DEA models: the fundamental role of interior points. *Journal of Productivity Analysis*, *38*, 121–130.

Krivonozhko, V. E., Førsund, F. R., & Lychev, A. V. (2014). Measurement of returns to scale using non-radial DEA models. *European Journal of Operational Research*, *232*, 664–670.

Lozano, S., & Villa, G. (2005). Determining a sequence of targets in DEA. *Journal of the operational research society*, *56*, 1439–1447.

Mehdiloozad, M., Mirdehghan, S. M., Sahoo, B. K., & Roshdi, I. (in press). On the identification of the global reference set in data envelopment analysis. *European Journal of Operational Research*.

Olesen, O. B., & Petersen, N. C. (1996). Indicators of ill-conditioned data sets and model misspecification in data envelopment analysis: An extended facet approach. *Management Science*, *42*, 205–219.

Olesen, O. B., & Petersen, N. C. (2003). Identification and use of efficient faces and facets in DEA. *Journal of productivity Analysis*, *20*, 323–360.

Pastor, J. T., & Aparicio, J. (2010). The relevance of DEA benchmarking information and the Least-Distance Measure: Comment. *Mathematical and Computer Modelling*. *52*, 397–399.

Portela, M. C. A. S., Borges, P. C., & Thanassoulis, E. (2003). Finding closest targets in non-oriented DEA models: The case of convex and non-convex technologies. *Journal of Productivity Analysis*. *19*, 251–269.

Roshdi, I., Mehdiloozad, M., & Margaritis, D. (2014). A linear programming based approach for determining maximal closest reference set in DEA. arXiv: 1407.2592 [math.OC].

Saaty T. (1980). *The Analytic Hierarchy Process*. New York: McGraw Hill.

Sueyoshi, T., & Sekitani, K. (2007a). The measurement of returns to scale under a simultaneous occurrence of multiple solutions in a reference set and a supporting hyperplane. *European Journal of Operational Research*, *181*, 549–570.